\documentclass[10pt]{amsart}
\usepackage{amssymb,amstext,amsmath,amscd,amsthm,amsfonts,enumerate,latexsym,stmaryrd,multicol,geometry,graphicx,mathrsfs,bm}
\usepackage[usenames]{color}
\usepackage[all]{xy}
\newtheorem{thm}{Theorem}[section]

\newtheorem{cor}[thm]{Corollary}
\theoremstyle{definition}

\newtheorem{rem}[thm]{Remark}

\theoremstyle{remark}

\newtheorem*{claim*}{Claim}
\newtheorem*{ac}{Acknowledgments}
\newtheorem*{conv}{Convention}

\numberwithin{equation}{thm}

\def\amp{\operatorname{amp}}

\def\db{\operatorname{\mathsf{D^b}}}
\def\depth{\operatorname{depth}}
\def\dim{\operatorname{dim}}
\def\ds{\operatorname{\mathsf{D_{sg}}}}

\def\edim{\operatorname{edim}}

\def\Ext{\operatorname{Ext}}

\def\ge{\geqslant}

\def\hom{\operatorname{Hom}}

\def\inf{\operatorname{inf}}

\def\le{\leqslant}

\def\m{\mathfrak{m}}

\def\perf{\operatorname{\mathsf{perf}}}

\def\sup{\operatorname{sup}}
\def\syz{\Omega}
\def\T{\mathcal{T}}

\def\X{\mathcal{X}}
\def\Y{\mathcal{Y}}

\def\ZZ{\mathbb{Z}}
\begin{document}
\title{Uniform dominance of Golod local rings}
\author{Toshinori Kobayashi}
\address[Kobayashi]{School of Science and Technology, Meiji University, 1-1-1 Higashimita, Tamaku, Kawasaki 214-8571, Japan}
\email{tkobayashi@meiji.ac.jp}
\author{Yuki Mifune}
\address[Mifune]{Graduate School of Mathematics, Nagoya University, Furocho, Chikusaku, Nagoya 464-8602, Japan}
\email{yuki.mifune.c9@math.nagoya-u.ac.jp}
\begin{abstract}
We prove that every Golod local ring $R$ is uniformly dominant. More precisely, we show that its dominant index is at most the embedding dimension of $R$. The proof is based on a numerical criterion for formality of bounded complexes whose homology is annihilated by the maximal ideal.
\end{abstract}
\maketitle
\section{Introduction}
The dominant index of a commutative noetherian local ring measures the number of extensions needed to construct its residue field from an arbitrary nonzero object of its singularity category, with shifts, finite direct sums, and direct summands allowed. 
Takahashi introduced uniform dominance and established it for several classes of local rings \cite{T2026}.
He had earlier asked whether every Golod local ring is dominant \cite[Question 9.4]{T2023}, and Kobayashi and Takahashi later posed the stronger question of whether every such ring is uniformly dominant \cite[Question 1.1(3)]{KT2026}.
The following theorem answers both questions.
\begin{thm}[Corollary \ref{cor_dom}]\label{mthm_int1}
Every Golod local ring $R$ is uniformly dominant with $\operatorname{dx}(R)\le\edim R$. 
\end{thm}
The proof of Theorem \ref{mthm_int1} is based on a numerical criterion for formality of bounded complexes. 
Let $(R,\m,k)$ be a local ring.
Denote by $\db(R)$ the bounded derived category of finitely generated $R$-modules.
For $C\in\db(R)$ and $n\in\ZZ$, set $\beta_n(C)=\dim_k\Ext_R^n(C,k)$ and $P_C(t)=\sum_{n\in\ZZ}\beta_n(C)t^n$. 
If $\m H_i(C)=0$ for all $i\in\ZZ$, we also set $h_i(C)=\dim_k H_i(C)$ for each $i\in\ZZ$ and $h_C(t)=\sum_{i\in\ZZ}h_i(C)t^i$.
The following theorem characterizes formality of such complexes in terms of their Poincar\'e series.
\begin{thm}[Theorem \ref{form}]\label{mthm_int2}
Let $(R,\m,k)$ be a local ring, and $C\in\db(R)$ a complex such that $\m H_i(C)=0$ for all $i\in\mathbb Z$. 
Then the following conditions are equivalent.
\begin{enumerate}[\rm(1)]
\item 
One has $P_C(t)=h_C(t)P_k(t)$.
\item 
The complex $C$ is formal, that is, $C\cong\bigoplus_{i\in\mathbb Z}H_i(C)[i]$ in $\db(R)$.
\end{enumerate}
\end{thm}
The key application of Theorem \ref{mthm_int2} is the formality of the Koszul complex of a Golod module over a Golod local ring; see Corollary \ref{cor_golod}. 
Combining this with Lescot's theorem, we prove Theorem \ref{mthm_int1}.
Theorem \ref{mthm_int1} also improves several known results by providing sharper bounds on dominant indices and establishing uniform dominance in greater generality; see Remark \ref{appl}.
\section{Formality and uniform dominance}
\begin{conv}
Throughout this paper, $(R,\m,k)$ denotes a commutative noetherian local ring. 
All subcategories are assumed to be strictly full. 
All complexes are chain complexes.
For a complex $C$ and an integer $n$, the shift $C[n]$ is defined by $(C[n])_i=C_{i-n}$ and $d_{C[n]}=(-1)^n d_C$. 
In particular, one has $H_i(C[n])=H_{i-n}(C)$.
\end{conv}
For a complex $C$ and an integer $n$, we write $\tau_{\ge n}C$ for the soft truncation of $C$ in degrees at least $n$ (see \cite[1.2]{CH2010}).
For a complex $C\in\db(R)$, set $\inf C=\inf\{i\in\ZZ\mid H_i(C)\ne0\}$ and $\sup C=\sup\{i\in\ZZ\mid H_i(C)\ne0\}$. 
For a nonzero complex $C\in\db(R)$, the {\em amplitude} of $C$ is defined by $\amp C=\sup C-\inf C$.
The symbol $\preccurlyeq$ denotes the coefficientwise inequality of Laurent series.
\begin{thm}\label{form}
Let $C\in\db(R)$ satisfy $\m H_i(C)=0$ for all $i\in\ZZ$. 
Then the following assertions hold.
\begin{enumerate}[\rm(1)]
\item
The inequality $P_C(t)\preccurlyeq h_C(t)P_k(t)$ holds.
\item 
The equality $P_C(t)=h_C(t)P_k(t)$ holds if and only if $C$ is formal in $\db(R)$.
\end{enumerate}
\end{thm}
\begin{proof}
We may assume that $C\ne0$. 
We prove the assertions by induction on $\amp C\ge0$. 
Put $a=\inf C$. 
If $\amp C=0$, then $C\cong H_a(C)[a]\cong k[a]^{\oplus h_a(C)}$, and both assertions follow.
Suppose that $\amp C>0$. 
Consider the canonical triangle
\[
\tau_{\ge a+1}C\to C\to H_a(C)[a]\xrightarrow{\delta}(\tau_{\ge a+1}C)[1].
\]
For each $n\in\ZZ$, let $\delta_n^\ast :\Ext_R^{n-1}(\tau_{\ge a+1}C,k)\to\Ext_R^n(H_a(C)[a],k)$ be the map induced by $\delta$.
Note that $\delta_n^\ast$ is naturally identified with $\hom_{\db(R)}(\delta,k[n])$.
Applying the cohomological functor $\hom_{\db(R)}(-,k)$ to the above triangle, we obtain the exact sequence
\[
\Ext_R^{n-1}(\tau_{\ge a+1}C,k)\xrightarrow{\delta_n^{\ast}}\Ext_R^n(H_a(C)[a],k)\to\Ext_R^n(C,k) \to\Ext_R^n(\tau_{\ge a+1}C,k)\xrightarrow{\delta_{n+1}^\ast}\Ext_R^{n+1}(H_a(C)[a],k)
\]
for each $n\in\ZZ$. 
Since $\Ext_R^n(H_a(C)[a],k)\cong\Ext_R^{n-a}(k,k)^{\oplus h_a(C)}$, exactness yields
\[
\beta_n(C)=\beta_n(\tau_{\ge a+1}C)+h_a(C)\beta_{n-a}(k)-\dim_k\operatorname{Im}\delta_n^\ast-\dim_k\operatorname{Im}\delta_{n+1}^{\ast}.
\]
Since $\amp(\tau_{\ge a+1}C)<\amp C$, the induction hypothesis gives $\beta_n(\tau_{\ge a+1}C)\le\sum_{i\neq a}h_i(C)\beta_{n-i}(k)$. 
Hence
\[
\beta_n(C)\le\beta_n(\tau_{\ge a+1}C)+h_a(C)\beta_{n-a}(k)\le\sum_{i\in\ZZ}h_i(C)\beta_{n-i}(k).
\]
This proves (1).

Suppose that $P_C(t)=h_C(t)P_k(t)$. 
Then the above inequalities are equalities for every $n\in\ZZ$. 
Hence $\delta_n^{\ast}=0$ for all $n\in\ZZ$, and $P_{\tau_{\ge a+1}C}(t)=h_{\tau_{\ge a+1}C}(t)P_k(t)$. 
By the induction hypothesis, $\tau_{\ge a+1}C$ is formal, so
\[
\tau_{\ge a+1}C\cong\bigoplus_{i\ne a}H_i(C)[i]\cong\bigoplus_{i\ne a}k[i]^{\oplus h_i(C)}
\]
in $\db(R)$.
Fix such an isomorphism, and let $\pi_{i,j}:(\tau_{\ge a+1}C)[1]\to k[i+1]^{\oplus h_i(C)}\to k[i+1]$ be the composition of projections for $i\ne a$ and $1\le j\le h_i(C)$. 
Since $\pi_{i,j}\in\hom_{\db(R)}((\tau_{\ge a+1}C)[1],k[i+1])\cong\Ext_R^{i}(\tau_{\ge a+1}C,k)$, we have $\pi_{i,j}\circ\delta=\delta_{i+1}^{\ast}(\pi_{i,j})=0$ for all $i,j$. Thus $\delta=0$, and the canonical triangle splits. Consequently, $C\cong\tau_{\ge a+1}C\oplus H_a(C)[a]\cong\bigoplus_{i\in\ZZ}H_i(C)[i]$, so $C$ is formal.

Conversely, suppose that $C$ is formal. Then $C\cong\bigoplus_{i\in\ZZ}k[i]^{\oplus h_i(C)}$ in $\db(R)$, and hence $\Ext_R^n(C,k)\cong\bigoplus_{i\in\ZZ}\Ext_R^{n-i}(k,k)^{\oplus h_i(C)}$. Thus $\beta_n(C)=\sum_{i\in\ZZ}h_i(C)\beta_{n-i}(k)$ for all $n\in\ZZ$, which gives $P_C(t)=h_C(t)P_k(t)$.
\end{proof}
\begin{rem}
When $C=K\otimes_R M$, where $K$ is the Koszul complex on a minimal generating set of $\m$ and $M$ is a nonzero finitely generated $R$-module, the inequality in Theorem \ref{form}(1) follows from \cite[Theorem 6.1.1]{AIM2006} applied to the identity map of $R$.
\end{rem}
Put $e=\edim R$, and let $K$ be the Koszul complex on a minimal generating set of $\m$. 
Set $D_R(t)=1-\sum_{i=1}^e\dim_k H_i(K)t^{i+1}$. 
Recall that $R$ is called a {\em Golod ring} if $P_k(t)=(1+t)^e/D_R(t)$.
A nonzero finitely generated $R$-module $M$ is called a {\em Golod module} if $P_M(t)=(\sum_{i=0}^e\dim_k H_i(K\otimes_R M)t^i)/D_R(t)$; see \cite[Section 5.3]{A1998}.
We apply Theorem \ref{form} to obtain the following characterization of Golod modules.
\begin{cor}\label{cor_golod}
Let $R$ be a Golod local ring, and let $K$ be the Koszul complex on a minimal generating set of $\m$. 
Then a nonzero finitely generated $R$-module $M$ is Golod if and only if $K\otimes_R M$ is formal in $\db(R)$.
\end{cor}
\begin{proof}
Let $P$ be a minimal free resolution of $M$. 
Since $K$ is a bounded complex of finite free modules, $P\otimes_R K$ is a minimal free model of $K\otimes_R M$. 
As $K_i\cong R^{\oplus\binom ei}$ for $0\le i\le e$, we have $\beta_n(K\otimes_R M)=\operatorname{rank}(P\otimes_R K)_n=\sum_{i=0}^e\binom ei\beta_{n-i}(M)$ for all $n\in\ZZ$.
Hence $P_{K\otimes_R M}(t)=(1+t)^eP_M(t)$.
Also, $\m H_i(K\otimes_R M)=0$ for all $i\in\ZZ$ by \cite[(1.5.3)]{AIM2006}.
By Theorem \ref{form} and the Golodness of $R$, the complex $K\otimes_R M$ is formal if and only if
\[
(1+t)^eP_M(t)=P_k(t)h_{K\otimes_R M}(t)=\frac{(1+t)^e}{D_R(t)}\sum_{i=0}^e\dim_k H_i(K\otimes_R M)t^i.
\]
Cancelling $(1+t)^e$, we obtain precisely the defining equality for $M$ to be a Golod module.
\end{proof}
To prove Theorem \ref{mthm_int1}, we recall the notion of generation in triangulated categories.
Let $\T$ be a triangulated category. 
For a subcategory $\X$ of $\T$, denote by $\langle\X\rangle_1^\T$ its closure under finite direct sums, direct summands and shifts. 
For subcategories $\X,\Y$ of $\T$, let $\X\ast\Y$ denote the subcategory consisting of objects $Z$ admitting an exact triangle $X\to Z\to Y\to X[1]$ with $X\in\X$ and $Y\in\Y$. 
For an object $X\in\T$, set $\langle X\rangle_0^\T=\{0\}$ and $\langle X\rangle_{n+1}^\T=\langle\langle X\rangle_n^\T\ast\langle X\rangle_1^\T\rangle_1^\T$ for $n\ge1$.

The {\em singularity category} of $R$ is defined as the Verdier quotient $\ds(R)=\db(R)/\perf(R)$, where $\perf(R)$ denotes the subcategory of perfect complexes.
Following \cite[Definition 3.4]{KT2026}, the {\em dominant index} of $R$ is defined by
\[
\operatorname{dx}(R)=\inf\{n\in\ZZ_{\ge-1}\mid k\in\langle X\rangle_{n+1}^{\ds(R)}\text{ for every }0\ne X\in\ds(R)\}.
\]
The local ring $R$ is called {\em uniformly dominant} if $\operatorname{dx}(R)<\infty$.
\begin{cor}\label{cor_dom}
Every Golod local ring $R$ is uniformly dominant with $\operatorname{dx}(R)\le\edim R$. 
In particular, Question 9.4 of \cite{T2023} and Question 1.1(3) of \cite{KT2026} have affirmative answers.
\end{cor}
\begin{proof}
Put $e=\edim R$, and let $0\ne X\in\ds(R)$. 
By \cite[Lemma 2.4(2)]{DT2015}, there exist a finitely generated $R$-module $M$ and an integer $m$ such that $X\cong\syz^e M[m]$ in $\ds(R)$. 
Lescot's theorem \cite[Theorem 5.3.2(2)]{A1998} implies that $\syz^e M$ is a nonzero Golod module.
By Corollary \ref{cor_golod}, the complex $K\otimes_R\syz^e M$ is formal in $\db(R)$. 
Since $H_0(K\otimes_R\syz^e M)\cong\syz^e M/\m\syz^e M\ne0$ by Nakayama's lemma, the residue field $k$ is a direct summand of $K\otimes_R\syz^e M$ in $\db(R)$.
On the other hand, the hard truncation yields $K\otimes_R\syz^e M\in\langle\syz^e M\rangle_{e+1}^{\db(R)}$. 
Passing to the singularity category, we obtain
\[
k\in\langle\syz^e M\rangle_{e+1}^{\ds(R)}=\langle M\rangle_{e+1}^{\ds(R)}=\langle X\rangle_{e+1}^{\ds(R)}.
\]
Thus $\operatorname{dx}(R)\le e$, proving the assertion.
\end{proof}

\begin{rem}\label{appl}
We record some consequences of Corollary \ref{cor_dom}
and compare them with known results.
Put $d=\dim R$.

\begin{enumerate}[\rm(1)]
\item
Suppose that the $\m$-adic completion $\widehat R$
is a hypersurface.
Then $R$ is Golod, and hence
$\operatorname{dx}(R)\le d+1$.
This extends \cite[Corollary 6.4]{KT2026}
by removing the completeness assumption,
while preserving the same upper bound.

\item
Suppose that
$\edim R-\depth R\le2$
and that $R$ is not a complete intersection.
Then $R$ is Golod by \cite[Propositions 5.2.5 and 5.3.4]{A1998}, so
\[
\operatorname{dx}(R)
\le\edim R
\le\depth R+2
\le d+2.
\]
In the Cohen--Macaulay case, this improves the
bound $6d+5$ in \cite[Theorem 6.14(2)]{KT2026}
to $d+2$. Moreover, neither the Cohen--Macaulay
assumption nor the infinitude of the residue
field is required.

\item
Suppose that $R$ is Cohen--Macaulay with
minimal multiplicity, that is,
$\mathrm e(R)=\edim R-d+1$.
Then $R$ is Golod, and consequently
$\operatorname{dx}(R)\le\edim R$.
For an infinite residue field, the sharper bound
$\operatorname{dx}(R)\le d$ was established in
\cite[Proposition 6.7(1b)]{KT2026}.
The present result establishes uniform dominance
without any restriction on the residue field,
although its numerical bound is weaker.
\item
Let
\[
R=k[\![x_1,\ldots,x_h]\!]/
I_2\begin{pmatrix}
f_{11}&\cdots&f_{1u}\\
f_{21}&\cdots&f_{2u}
\end{pmatrix},
\]
where $k$ is a field and
$f_{ij}\in(x_1,\ldots,x_h)$.
Suppose that $\dim R=h-u+1$.
Then $R$ is Golod, and hence
$\operatorname{dx}(R)\le h$.
Indeed, the generic determinantal ring defined
by the $2$-minors of a $2\times u$ matrix has
minimal multiplicity and is therefore Golod.
The ring $R$ is obtained from a formal power
series extension of this ring by a regular
sequence whose elements are linearly independent
modulo the square of the maximal ideal.
Thus $R$ is Golod by
\cite[Proposition 5.2.4(3)]{A1998}.
This substantially improves the bound in
\cite[Corollary 6.8]{T2026}.
\item
Let $S$ be a regular ring, let
$f_1,\ldots,f_n$ be an $S$-regular sequence,
and let $s\ge 2$.
Set $R=S/(f_1,\ldots,f_n)^s$.
Since powers of complete intersection ideals
in regular local rings define Golod rings
by \cite[Corollary 3.3]{GS1980},
every localization $R_{\mathfrak p}$ is Golod.
Consequently,
$\operatorname{dx}(R_{\mathfrak p})
\le \edim R_{\mathfrak p}$ for every $\mathfrak p\in\operatorname{Spec} R$.
This strengthens \cite[Proposition 5.1]{KMOT2026}
from local dominance to uniform dominance
at every localization.

\end{enumerate}
\end{rem}
\begin{ac}
Kobayashi was partly supported by JSPS Grant-in-Aid for Early-Career Scientists 25K17240.
Mifune was partly supported by Grant-in-Aid for JSPS Fellows 25KJ1386.
OpenAI's ChatGPT was used for mathematical discussions during the preparation of this paper.
\end{ac}


\begin{thebibliography}{99}
\bibitem{A1998}
{\sc L. L. Avramov}, Infinite free resolutions, {\em Six lectures on commutative algebra (Bellaterra, 1996)}, 1--118, Progr. Math. {\bf 166}, {\em Birkh\"auser, Basel}, 1998.
\bibitem{AIM2006}
{\sc L. L. Avramov; S. B. Iyengar; C. Miller}, Homology over local homomorphisms, {\em Amer. J. Math.} {\bf 128} (2006), no. 1, 23--90.
\bibitem{CH2010}
{\sc L. W. Christensen; H. Holm}, Algebras that satisfy Auslander's condition on vanishing of cohomology, {\em Math. Z.} {\bf 265} (2010), no. 1, 21--40.
\bibitem{DT2015}
{\sc H. Dao; R. Takahashi}, Upper bounds for dimensions of singularity categories, {\em C. R. Math. Acad. Sci. Paris} {\bf 353} (2015), no. 4, 297--301.
\bibitem{GS1980}
{\sc E. Gover; P. Salmon}, On a class of Golod homomorphisms. {\em Math. Scand.} {\bf 46} (1980), no. 1, 5--14.
\bibitem{KMOT2026}
{\sc K. Kimura; Y. Mifune; Y. Otake; R. Takahashi}, On strongly G-regular rings, preprint (2026), {\tt arXiv:2608.08228v1}.
\bibitem{KT2026}
{\sc T. Kobayashi; R. Takahashi}, On the ubiquity of uniformly dominant local rings, preprint (2026), {\tt arXiv:2603.10810v1}.
\bibitem{T2023}
{\sc R. Takahashi}, Dominant local rings and subcategory classification, {\em Int. Math. Res. Not. IMRN} (2023), no. 9, 7259--7318.
\bibitem{T2026}
{\sc R. Takahashi}, Uniformly dominant local rings and Orlov spectra of singularity categories, {\em Math. Z.} {\bf 312} (2026), no. 4, Paper No. 117, 27 pp.
\end{thebibliography}
\end{document}